# BOREL THEOREMS FOR RANDOM MATRICES FROM THE CLASSICAL COMPACT SYMMETRIC SPACES[1]

By Benoît Collins and Michael Stolz

*Université Claude Bernard Lyon 1 and Ruhr-Universität Bochum*

We study random vectors of the form $(\mathrm{Tr}(A^{(1)}V), \ldots, \mathrm{Tr}(A^{(r)}V))$, where $V$ is a uniformly distributed element of a matrix version of a classical compact symmetric space, and the $A^{(\nu)}$ are deterministic parameter matrices. We show that for increasing matrix sizes these random vectors converge to a joint Gaussian limit, and compute its covariances. This generalizes previous work of Diaconis et al. for Haar distributed matrices from the classical compact groups. The proof uses integration formulas, due to Collins and Śniady, for polynomial functions on the classical compact groups.

## 1. Introduction.
One way to extend the scope of random matrix theory is to revisit its roots in physics. Classically, random matrices were considered by physicists as toy models for operators associated to complicated quantum systems. The matrix ensemble was supposed to be compatible with the symmetries of the quantum system. Analyzing the action of symmetry groups, Freeman Dyson [10] argued in 1962 that the good objects to study were the Gaussian Orthogonal, Unitary and Symplectic Ensembles (GOE, GUE, GSE) as models for Hamiltonians, and the Circular Orthogonal, Unitary and Symplectic Ensembles (COE, CUE, CSE) as models for time evolution operators. The geometric fact behind Dyson's result is that the matrix ensembles in question are associated to classical symmetric spaces.

For concreteness, let us review the case of GOE/COE: GOE($n$) is a probability measure on the space of real symmetric $n \times n$ matrices such that the matrix entries on and above the diagonal are independent centered Gaussian random variables. What is more, it is assumed that the variance of the diagonal entries is twice the variance of the off-diagonal entries. This guarantees that the Lebesgue density of GOE($n$) factorizes with respect to the

Received December 2006; revised December 2006.
[1]Supported by Deutsche Forschungsgemeinschaft via SFB/TR 12.
AMS 2000 subject classifications. Primary 15A52, 60F05; secondary 60B15, 43A75.
Key words and phrases. Random matrices, symmetric spaces, central limit theorem, matrix integrals, classical invariant theory.

trace and is hence invariant under conjugation by orthogonal matrices. Let us give an abstract account of these well-known facts: Write

$$\mathfrak{u} := \mathfrak{u}_n := \{X \in \mathbb{C}^{n \times n} \colon \overline{X}' = -X\},$$

the real vector space of skew-Hermitian complex $n \times n$ matrices. $\mathfrak{u}$ admits a real linear involution $\theta(X) = \overline{X}$. Then we have the eigenspace decomposition

$$\mathfrak{u} = \mathfrak{k} \oplus \mathfrak{p},$$

where

$$\mathfrak{k} = \{X \in \mathbb{R}^{n \times n} \colon X' = -X\}$$

is the $+1$ eigenspace and

$$\mathfrak{p} = \{X \in i\mathbb{R}^{n \times n} \colon X' = X\}$$

is the $-1$ eigenspace. Note that

$$i\mathfrak{p} = \{X \in \mathbb{R}^{n \times n} \colon X' = X\}$$

is the space on which $\mathrm{GOE}(n)$ is supported. Applying the matrix exponential to the elements of $\mathfrak{u}$, we have

$$e^X \overline{e^{X'}} = e^X e^{\overline{X}'} = e^X e^{-X} = I,$$

hence $e^X$ is unitary. In fact, any unitary matrix arises in this way, and $\mathfrak{u}_n$ is the Lie algebra of the unitary group $\mathrm{U}_n$. $\mathfrak{k}$ exponentiates to the special orthogonal group $\mathrm{SO}_n$ and is its Lie algebra. So, on this abstract level, $\mathrm{GOE}(n)$ is a probability distribution on $i\mathfrak{p}$, which is invariant under the action of the Lie group corresponding to $\mathfrak{k}$.

For $X \in \mathfrak{p}$, the image $e^X$ is both unitary and symmetric, and it turns out that

$$\exp(\mathfrak{p}) = \{gg' \colon g \in \mathrm{U}_n\} =: \mathcal{S}.$$

For $g \in \mathrm{U}_n$ write $\Phi(g) = gg'$. Then $\Phi^{-1}(1) = \mathrm{O}_n$. We can thus view $\Phi$ as a bijection of $\mathrm{U}_n/\mathrm{O}_n$ onto $\mathcal{S}$. The uniform distribution on $\mathrm{U}_n/\mathrm{O}_n$ is Dyson's circular ensemble COE. Note that $\mathrm{O}_n \subset \mathrm{U}_n$ is the set of fixed points of the involution $\tilde{\theta}(g) = (g^{-1})'$ of $\mathrm{U}_n$. Coset spaces of this kind arise when Riemannian symmetric spaces (RSS) are represented as homogeneous spaces of Lie groups, and conversely, modulo some technicalities, coset spaces of this kind are RSS (see [14] for details). In this sense we can view $\mathcal{S}$ and $\mathfrak{p}$ (hence $i\mathfrak{p}$) as representing an RSS.

As a matter of fact, the matrix ensembles from Dyson's list are of Cartan class A, AI, AII; that is, $\mathrm{GUE}(n)$ lives on $i\mathfrak{u}$, which is nothing else than the space of hermitian $n \times n$ matrices, $\mathrm{GOE}(n)$ is as above, and $\mathrm{GSE}(n)$ lives on the $i\mathfrak{p}$ part in $\mathfrak{u}_{2n}$, this time with respect to the involution $\theta(X) = J_n \overline{X} J_n^{-1}$, where $J_n$ is as in (1) below. In recent years, it has emerged that the framework that led to Dyson's list is not flexible enough to meet the needs of condensed matter physics; see, for example, [1]. A classification program, in the spirit of Dyson, has been carried out by Heinzner, Huckleberry and Zirnbauer in [13], to the effect that all spaces $i\mathfrak{g}$, $i\mathfrak{p}$ arise naturally as spaces of (toy models for) Hamiltonian operators in mesoscopic physics, where $\mathfrak{g}$ is any classical compact Lie algebra and $\mathfrak{p}$ is the $-1$ eigenspace of a classical "Cartan involution" of $\mathfrak{g}$ (see Section 2 below). By the same token, all classical compact Lie groups and spaces $\mathcal{S}$ are needed to model time evolution operators.

This paper is devoted to the $\mathcal{S}$ side of the story. For the $i\mathfrak{g}, i\mathfrak{p}$ side see [11]; for other aspects of the $\mathcal{S}$ side see [9]. We endow $\mathcal{S}$ with the uniform distribution and investigate the behavior of the matrix entries of a random element of $\mathcal{S}$ as matrix size tends to infinity. We obtain Gaussian limits for wide classes of vectors of linear combinations of matrix entries. The archetypical result in this direction is due to Émile Borel, who proved a century ago that for a uniformly distributed point $X = (X_1, \ldots, X_n)$ on the unit sphere $S^{n-1} \subset \mathbb{R}^n$, the scaled first coordinate $\sqrt{n}X_1$ converges weakly to the standard normal distribution as $n \to \infty$. Further-reaching results about the entries of matrices from the classical groups have been obtained in [5, 6] and [15], the first reference containing illuminating remarks on the history of the topic. To prove the weak convergence results of this paper, we invoke the method of moments. To make it operational, we use what we term "extended Wick calculus." The key ideas can be traced back to Weingarten [17]. They have been developed by Collins [3] and Collins and Śniady [4]. There is a certain family resemblance between this calculus and the Diaconis–Evans–Mallows–Ram–Shahshahani approach to random matrix eigenvalues [7, 8], in that both exploit in different ways the invariant theory of the classical groups (see [16] and Section 3 below).

The paper is organized as follows: In Section 2, we explain the basic setup and notation and state the main results. Then, in Section 3, we review extended Wick calculus and explain its connection with the Diaconis–Shahshahani integral. Finally, the main results are proven in Section 4.

## 2. Notation and statement of the main results.
For $n \in \mathbb{N}$ let $G = G_n$ be one of the unitary group $U_n$, the (compact) orthogonal group $O_n$ and the (compact) symplectic group $Sp_{2n}$. $g'$ denotes the transpose of a matrix $g$, and $g^*$ denotes conjugate transpose. For $n \in \mathbb{N}$, write $I_n \in \mathbb{C}^{n \times n}$ for the identity matrix. For $p, q, n \in \mathbb{N}, p + q = n$, set

$$(1) \qquad I_{p,q} = \begin{pmatrix} I_p & 0 \\ 0 & -I_q \end{pmatrix}, \qquad J = J_n = \begin{pmatrix} 0 & -I_n \\ I_n & 0 \end{pmatrix}$$

and
$$\mathcal{I}_{p,q} = \begin{pmatrix} \mathrm{I}_{p,q} & 0 \\ 0 & \mathrm{I}_{p,q} \end{pmatrix}.$$

Note that $\mathrm{I}_{p,q} = \mathrm{I}_{p,q}^{-1}$ and $\mathrm{J}' = \mathrm{J}^{-1}$. We consider pairs $(G_n, \theta_{\mathcal{C}})$, where $\theta_{\mathcal{C}}$ is an involutive automorphism of $G_n$, given as follows:

- For $G_n = \mathrm{U}_{\delta n}$, we set
$$\theta_{\mathrm{AI}}(g) = \overline{g},$$
$$\theta_{\mathrm{AII}}(g) = J'\overline{g}J,$$
$$\theta_{\mathrm{AIII}}(g) = \mathrm{I}_{p,q}g\mathrm{I}_{p,q}.$$

- For $G_n = \mathrm{O}_{\delta n}$, we set
$$\theta_{\mathrm{BDI}}(g) = \mathrm{I}_{p,q}g\mathrm{I}_{p,q},$$
$$\theta_{\mathrm{DIII}}(g) = J'gJ.$$

- For $G_n = \mathrm{Sp}_{2n}$, we set
$$\theta_{\mathrm{CI}}(g) = \mathrm{I}_{n,n}g\mathrm{I}_{n,n},$$
$$\theta_{\mathrm{CII}}(g) = \mathcal{I}_{p,q}g\mathcal{I}_{p,q}.$$

Here, $\delta = \delta_{\mathcal{C}} = 2$ for $\mathcal{C} = \mathrm{AII}, \mathrm{DIII}$ and $\delta = 1$ otherwise. Write $K = K_n = K_n^{\mathcal{C}}$ for the fixed subgroups of the automorphisms $\theta_{\mathcal{C}}$. Concretely, $K_n^{\mathrm{AI}} = \mathrm{O}_n$, $K_n^{\mathrm{AII}} = \mathrm{Sp}_{2n}$, $K_n^{\mathrm{AIII}} = \mathrm{U}_p \times \mathrm{U}_q$, $K_n^{\mathrm{BDI}} = \mathrm{O}_p \times \mathrm{O}_q$, $K_n^{\mathrm{DIII}} = \mathrm{U}_n \subset \mathrm{O}_{2n}$, $K_n^{\mathrm{CI}} = \mathrm{U}_n \subset \mathrm{Sp}_{2n}$, $K_n^{\mathrm{CII}} = \mathrm{Sp}_{2p} \times \mathrm{Sp}_{2q}$.

Then the list of homogeneous spaces $G_n/K_n^{\mathcal{C}}$ is *cum grano salis* the list of classical compact symmetric spaces, as classified by Cartan. But it should be emphasized that we follow standard practice in probability theory in that we consider the full groups $\mathrm{U}_n$ and $\mathrm{O}_n$ rather than $\mathrm{SU}_n$ and $\mathrm{SO}_n$, which would be the natural objects to consider in Lie theory. We will refer to $(G_n, \theta_{\mathcal{C}})$ as symmetry class $\mathcal{C}$, and we will say that $\mathrm{U}_n$ is of class A, $\mathrm{O}_n$ is of class B/D and $\mathrm{Sp}_n$ is of class C. The differentials of the $\theta_{\mathcal{C}}$ are the involutive automorphisms of the Lie algebra $\mathfrak{g}_n$ of $G_n$ which were discussed in the [Introduction](). Write

$$\Phi_{\mathcal{C}}(g) = \Phi_{\mathcal{C}}^{(n)}(g) = g(\theta_{\mathcal{C}}(g))^{-1} \quad \text{and} \quad \mathcal{S}_{\mathcal{C}} = \mathcal{S}_{\mathcal{C}}^{(n)} = \Phi_{\mathcal{C}}(G_n) \subset G_n.$$

We will refer to $\Phi_{\mathcal{C}}$ as Cartan embedding and to $\mathcal{S}_{\mathcal{C}}$ as the Cartan-embedded version of the compact symmetric space $G_n/K_n^{\mathcal{C}}$ of class $\mathcal{C}$. To make the connection with the preliminary discussion of GOE/COE in the [Introduction](), note that the tangent space in 1 to $\mathcal{S}_{\mathcal{C}}$ is the $(-1)$ eigenspace $\mathfrak{p}$ of the differential of $\theta_{\mathcal{C}}$. Adopting physics terminology, we call A, AI, AII the "Wigner–Dyson classes," AIII, BDI and CII the "chiral classes" and all others "superconductor" or "Bogolioubov–de Gennes (BdG)" classes. Observe that the chiral classes are special in that they depend not only on the overall matrix size $n$ (resp. $2n$), but also on parameters $p = p(n)$, $q = q(n)$ with $p + q = n$.

For completeness, define $\Phi_{\mathrm{A}}^{(n)}, \Phi_{\mathrm{B/D}}^{(n)}, \Phi_{\mathrm{C}}^{(n)}$ as identity on $\mathrm{U}_n, \mathrm{O}_n, \mathrm{Sp}_{2n}$, respectively. Write $\lambda_G = \lambda_{G_n}$ for Haar measure on $G_n$, normalized as a probability measure. Write $V = V_{\mathcal{C}} = V_{\mathcal{C}}^{(n)}$ for a random element of $\mathcal{S}_{\mathcal{C}}$, whose distribution is the push-forward of $\lambda_G$ under $\Phi_{\mathcal{C}}$. It is convenient to call classes A, AI, AII, AIII the "A classes," B/D, BDI, DIII the "B/D classes," and C, CI, CII the "C classes." Set $\mathbb{K} = \mathbb{C}$ for the A-classes, $\mathbb{K} = \mathbb{R}$ for the B/D-classes and $\mathbb{K} = \mathbb{H}$ for the C classes, $\mathbb{H}$ denoting the skew field of quaternions. We embed the space $\mathbb{H}^{n \times n}$ of quaternionic matrices into $\mathbb{C}^{2n \times 2n}$ as

$$\mathbb{H}^{n \times n} = \left\{ \begin{pmatrix} X & Y \\ -\overline{Y} & \overline{X} \end{pmatrix} : X, Y \in \mathbb{C}^{n \times n} \right\}.$$

Recall that $\mathrm{Sp}_{2n} = \mathbb{H}^{n \times n} \cap \mathrm{U}_{2n}$ and $\mathrm{O}_n = \mathbb{R}^{n \times n} \cap \mathrm{U}_n$. We endow $\mathbb{K}^{n \times n}$ with the scalar product $(A, B) \mapsto \mathrm{Tr}(AB^*)$. Write $W_{\mathcal{C}} = W_{\mathcal{C}}^{(n)}$ for the smallest $*$-closed subspace of $\mathbb{K}^{\delta n \times \delta n}$ which contains $\mathcal{S}_{\mathcal{C}}^{(n)}$, and $P_{\mathcal{C}}^{(n)}$ for the orthogonal projection of $\mathbb{K}^{\delta n \times \delta n}$ onto $W_{\mathcal{C}}^{(n)}$.

THEOREM 1. *Let $\mathcal{C}$ be a Wigner–Dyson or superconductor/BdG symmetry class. For $\nu = 1, \ldots, r$ and $n \in \mathbb{N}$ let $A^{(\nu,n)} \in \mathbb{K}^{\delta n \times \delta n}$ be a fixed deterministic matrix such that for all $\nu$ there exists a positive limit*

$$(2) \qquad \lim_{n \to \infty} \frac{1}{n} \mathrm{Tr}[P_{\mathcal{C}}^{(n)}(A^{(\nu,n)})(P_{\mathcal{C}}^{(n)}(A^{(\nu,n)}))^*].$$

*Then*

(i) *As $n \to \infty$, $\mathrm{Re}(\mathrm{Tr}(A^{(\nu,n)} V_{\mathcal{C}}^{(n)}))$ converges weakly to a random variable $\mathcal{T}_\nu$ ($\nu = 1, \ldots, r$).*

(ii) *$(\mathcal{T}_1, \ldots, \mathcal{T}_r)$ has a centered joint Gaussian distribution, and there exists a constant $\gamma_{\mathcal{C}} > 0$ such that for all $\mu, \nu = 1, \ldots, r$*

$$(3) \qquad \mathrm{Cov}(\mathcal{T}_\mu, \mathcal{T}_\nu) = \lim_{n \to \infty} \frac{\gamma_{\mathcal{C}}}{n} \mathrm{Re}\{\mathrm{Tr}[P_{\mathcal{C}}^{(n)}(A^{(\mu,n)})(P_{\mathcal{C}}^{(n)}(A^{(\nu,n)}))^*]\}.$$

THEOREM 2. *Let $\mathcal{C}$ be a chiral symmetry class. For $\nu = 1, \ldots, r$ and $n \in \mathbb{N}$ let $A^{(\nu,n)} \in \mathbb{K}^{n \times n}$ be a fixed deterministic matrix such that for all $\nu$ there exists a positive limit (2). Assume furthermore that the sequences $(p(n))_n$, $(q(n))_n$ [$n = p(n) + q(n)$] have the property that for all $\nu$ there exists the limit*

$$t_{\mathcal{C},\nu} := \lim_{n \to \infty} \frac{p(n) - q(n)}{n} \mathrm{Tr}(P_{\mathcal{C}}^{(n)}(A^{(\nu,n)}) I_{p(n),q(n)}),$$

*if $\mathcal{C}$ = BDI and AIII, and*

$$t_{\mathrm{CII},\nu} := \lim_{n \to \infty} 2\frac{p(n) - q(n)}{n} \operatorname{Tr}(P_{\mathrm{CII}}^{(n)}(A^{(\nu,n)})\mathcal{I}_{p(n),q(n)}),$$

*if $\mathcal{C}$ = CII. Then*

(i) *As $n \to \infty$, $\operatorname{Re}(\operatorname{Tr}(A^{(\nu,n)}V_{\mathcal{C}}^{(n,p(n),q(n))}))$ converges weakly to a random variable $\mathcal{T}_\nu$ ($\nu = 1, \ldots, r$).*

(ii) *$(\mathcal{T}_1, \ldots, \mathcal{T}_r)$ has a Gaussian distribution with mean vector $(t_{\mathcal{C},1}, \ldots, t_{\mathcal{C},r})$ and covariance matrix given by (3) for a suitable $\gamma_{\mathcal{C}} > 0$.*

REMARK 2.1.

(i)  $\gamma_{\mathcal{C}} = \frac{1}{2}$ for $\mathcal{C}$ = A, AII,
   $\gamma_{\mathcal{C}} = 1$ for $\mathcal{C}$ = AI, AIII, B/D, DIII, C,
   $\gamma_{\mathcal{C}} = 2$ for $\mathcal{C}$ = BDI,
   $\gamma_{\mathcal{C}} = 4$ for $\mathcal{C}$ = CI, CII.

(ii) For the B/D and C classes, $\operatorname{Tr}(A^{(\nu,n)}V_{\mathcal{C}}^{(n)})$ is in fact real. It follows immediately from the explicit description of the spaces $W_{\mathcal{C}}$ in Section 4 that the $t_{\mathcal{C},\nu}$ are real.

We close this section with some technical notation that will be important for all subsequent sections. For sets $\Omega$, $\Omega'$ write $\mathcal{F}(\Omega, \Omega') := \{\varphi : \Omega \to \Omega'\}$. If $\Omega = \{1, 2, \ldots, k\}$, $\Omega' = \{1, 2, \ldots, n\}$, write $\mathcal{F}(k, n) := \mathcal{F}(\Omega, \Omega')$. If $\#\Omega$ is even, write $\mathcal{M}(\Omega)$ for the set of all pair partitions of $\Omega$, that is, partitions of $\Omega$ all of whose blocks consist of exactly two elements. If $\Omega = \{1, 2, \ldots, 2k\}$, write $\mathcal{M}(2k) := \mathcal{M}(\Omega)$. For $\mathfrak{m} \in \mathcal{M}(\Omega)$ write $\mathcal{F}(\mathfrak{m}, \Omega')$ for the set of all $\varphi \in \mathcal{F}(\Omega, \Omega')$ which are constant on the blocks of $\mathfrak{m}$. We sometimes write functions in the right operator fashion, that is, $j\varphi$ in the place of $\varphi(j)$.

3. **Extended Wick calculus.** We will prove Theorems 1 and 2 via convergence of cumulants, hence of moments. To make this strategy operational, we need tools for the asymptotic evaluation of integrals of polynomial functions over the classical compact groups. These are available from recent work of Collins and Śniady [4]. To give a self-contained statement of these integral formulas, and to insert a sign factor into the symplectic integration formula which is missing in [4], we have to explain how these formulas are rooted in the invariant theory of the classical groups. As an extra benefit, this makes it possible to shed some light on how our tools relate to the Diaconis–Shahshahani integral. Let us start with this point, and focus on the conceptually simplest case, the orthogonal one. Consider a random variable $\Gamma_n$ with values in $\mathrm{O}_n$, distributed according to Haar measure. In order to prove via the method of moments that the random vector

$$(\operatorname{Tr}(\Gamma_n), \operatorname{Tr}(\Gamma_n^2), \ldots, \operatorname{Tr}(\Gamma_n^r))$$

converges to a Gaussian limit as $n \to \infty$, one has to evaluate matrix integrals of the form

$$(4) \qquad \int (\mathrm{Tr}(g))^{a_1} (\mathrm{Tr}(g^2))^{a_2} \cdots (\mathrm{Tr}(g^r))^{a_r} \lambda_{\mathrm{O}_n}(dg)$$

$[r \in \mathbb{N}, a = (a_1, \ldots, a_r) \in \mathbb{N}_0^r]$. As it turns out, the integral (4) has an interpretation in terms of the invariant theory of the group $\mathrm{O}_n$. Specifically, let $V$ be an $n$-dimensional complex vector space with scalar product $\langle \cdot, \cdot \rangle$ and symmetric bilinear form $\beta$, such that $\mathrm{O}(n, \mathbb{C})$ is the group of isometries of $\beta$ and $\mathrm{O}_n = \mathrm{O}(n, \mathbb{C}) \cap \mathrm{U}_n$, where $\mathrm{U}_n$ is the group of unitary matrices w.r.t. $\langle \cdot, \cdot \rangle$. Write $\rho$ for the natural representation of $\mathrm{O}(n, \mathbb{C})$ (or $\mathrm{O}_n$) on $V$. On the tensor power $V^{\otimes k}$ consider the representation $\rho_k = \rho^{\otimes k}$, which is given by

$$\left( \bigotimes_{j=1}^{k} v_j \right) \rho_k(g) := \bigotimes_{j=1}^{k} v_j \rho(g).$$

Define a representation $\sigma_k$ of the symmetric group $\mathrm{S}_k$ on $V^{\otimes k}$ via

$$\left( \bigotimes_{j=1}^{k} v_j \right) \sigma_k(s) := \bigotimes_{j=1}^{k} v_{js^{-1}}.$$

Then the integral (4) equals

$$(5) \qquad \mathrm{Tr}(\sigma_k(s)|_{[V^{\otimes k}]^{\mathrm{O}_n}}),$$

where $[V^{\otimes k}]^{\mathrm{O}_n} = \{ v \in V^{\otimes k} : v \rho_k(g) = v \; \forall g \in \mathrm{O}_n \}$ is the space of invariants in $V^{\otimes k}$ w.r.t. the representation $\rho_k$, and $s \in \mathrm{S}_k$ has cycle type $(1^{a_1} 2^{a_2} \cdots r^{a_r})$, whence $k = \sum_{j=1}^{r} j a_j$. See [16] for details.

Now we turn to the problem whose solution in [4] is the basic tool for our present purposes. Given $\varphi, \psi \in \mathcal{F}(k, n)$ (see notation at the end of Section 2), what can one say about the integral

$$(6) \qquad \int \prod_{j=1}^{k} g_{j\varphi, j\psi} \lambda_{\mathrm{O}_n}(dg).$$

If $e_j$ $(j = 1, \ldots, n)$ is an orthonormal basis of $V$, $\varphi \in \mathcal{F}(k, n)$, write $e_\varphi := \bigotimes_{j=1}^{k} e_{j\varphi}$. Then $\{ e_\varphi : \varphi \in \mathcal{F}(k, n) \}$ is a basis of $V^{\otimes k}$, orthonormal w.r.t. the extension of $\beta$ given by

$$(7) \qquad \beta \left( \bigotimes_{j=1}^{k} v_j, \bigotimes_{j=1}^{k} w_j \right) = \prod_{j=1}^{k} \beta(v_j, w_j).$$

Observe that (6) equals

$$\int \prod_{j=1}^{k} \beta(e_{j\varphi} g, e_{j\psi}) \, dg = \int \beta(e_\varphi \rho_k(g), e_\psi) \, dg = \beta \left( \int e_\varphi \rho_k(g) \, dg, e_\psi \right).$$

Now, $v \mapsto \int v \rho_k(g) \, dg$ is the orthogonal projection of $V^{\otimes k}$ onto $[V^{\otimes k}]^{\mathrm{O}_n}$, hence the integral (6) is the $e_\psi$-coordinate of the projection of $e_\varphi$ onto the space of invariants. This clarifies in algebraic terms the family resemblance between the matrix integrals (4) and (6). We will pursue this idea a bit further and show that all ingredients of the orthogonal integration formula from [4] can be given transparent definitions in terms of the function spaces $\mathcal{F}(\mathfrak{m}, n)$. The link to the approach in [4], via Brauer algebras, is provided by the fact that if the underlying vector space has sufficiently large dimension, the $\mathcal{F}(\mathfrak{m}, n)$ describe a basis of the space of tensor invariants of the orthogonal group. Setting $\mathfrak{m}_0 = \{\{1, 2\}, \ldots, \{2l - 1, 2l\}\} \in \mathcal{M}(2l)$,

$$\theta_l = \sum_{\varphi \in \mathcal{F}(\mathfrak{m}_0, n)} e_\varphi$$

and

$$\theta_l^{\mathrm{S}_{2l}} = \{\theta_l s : s \in \mathrm{S}_{2l}\},$$

Weyl's First Fundamental Theorem of the invariant theory of the orthogonal group says that for $k$ odd, $[V^{\otimes k}]^{\mathrm{O}_n}$ is trivial, whereas for $k = 2l$ even, it is spanned by the $\mathrm{S}_{2l}$-orbit $\theta_l^{\mathrm{S}_{2l}}$. In fact, for $n \geq l$, this orbit is a basis of the space of invariants, and it is parametrized by $\mathcal{M}(2l)$ ([16], Corollary 3.11, Lemma 3.12). Specifically, for $\mathfrak{m} \in \mathcal{M}(2l)$, there exists $s_{\mathfrak{m}} \in \mathrm{S}_{2l}$ such that

$$\theta_l s_{\mathfrak{m}} = \sum_{\varphi \in \mathcal{F}(\mathfrak{m}, n)} e_\varphi.$$

Set

$$c(\mathfrak{m}, \mathfrak{n}) := \beta(\theta_l s_{\mathfrak{m}}, \theta_l s_{\mathfrak{n}}).$$

For $n \geq l$ the matrix $C := (c(\mathfrak{m}, \mathfrak{n}))_{\mathfrak{m}, \mathfrak{n} \in \mathcal{M}(2l)}$ is invertible, and its inverse

$$(8) \qquad (\mathrm{Wg}_{\mathrm{O}_n}(\mathfrak{m}, \mathfrak{n}))_{\mathfrak{m}, \mathfrak{n}} := C^{-1}$$

coincides with the orthogonal Weingarten function as defined in [4] (where the reader will also find a definition for the case $n < l$). To see this, note that if $T$ is the $\mathrm{O}_n$-equivariant isomorphism from $V^{\otimes 2l}$ to $\mathrm{End}(V^{\otimes l})$ given by $u(T(v_1 \otimes \cdots \otimes v_{2l})) = \beta(u, v_2 \otimes v_4 \otimes \cdots \otimes v_{2l}) v_1 \otimes v_3 \otimes \cdots \otimes v_{2l-1}$ (see [12], Section 4.3.2), then $\rho_B(\mathfrak{p})$, as defined in [4], equation (15), equals $(\theta_l s_{\mathfrak{m}}) T$, where the Brauer diagram $\mathfrak{p}$ equals the pair partition $\mathfrak{m}$ if the upper and lower rows of the diagram are identified with even and odd numbers, respectively.

Now, $\int e_\varphi \rho_k(g) \, dg$ being an invariant, we see that it vanishes if $k$ is odd, and that for $k = 2l$ it can be written as $\sum_{\mathfrak{m} \in \mathcal{M}(k)} \alpha_{\mathfrak{m}}(\theta_l s_{\mathfrak{m}})$. Then we have

$$\alpha_{\mathfrak{m}} = \sum_{\mathfrak{n} \in \mathcal{M}(2l)} \beta(e_\varphi, \theta_l s_{\mathfrak{n}}) \mathrm{Wg}_{\mathrm{O}_n}(\mathfrak{m}, \mathfrak{n}),$$

hence

$$\beta\Big(\int e_\varphi \rho_k(g)\,dg, e_\psi\Big) = \sum_{\mathfrak{m}\in\mathcal{M}(2l)} \alpha_{\mathfrak{m}}\beta(\theta_l s_{\mathfrak{m}}, e_\psi)$$

$$= \sum_{\mathfrak{m},\mathfrak{n}} \beta(e_\varphi, \theta_l s_{\mathfrak{n}}) \mathrm{Wg}_{\mathrm{O}_n}(\mathfrak{m},\mathfrak{n})\beta(\theta_l s_{\mathfrak{m}}, e_\psi)$$

$$= \sum_{\mathfrak{m},\mathfrak{n}} \mathrm{Wg}_{\mathrm{O}_n}(\mathfrak{m},\mathfrak{n}) 1_{\mathcal{F}(\mathfrak{m},n)}(\varphi) 1_{\mathcal{F}(\mathfrak{n},n)}(\psi).$$

So we obtain the following

PROPOSITION 3.1. *Let $\varphi,\ \psi\in\mathcal{F}(k,n)$. Then*

$$\int_{\mathrm{O}_n} \prod_{j=1}^k g_{j\varphi,j\psi}\lambda_{\mathrm{O}_n}(dg) = \sum_{\mathfrak{m},\mathfrak{n}\in\mathcal{M}(k)} 1_{\mathcal{F}(\mathfrak{m},n)}(\varphi) 1_{\mathcal{F}(\mathfrak{n},n)}(\psi)\mathrm{Wg}_{\mathrm{O}_n}(\mathfrak{m},\mathfrak{n}).$$

*In particular, the integral vanishes if $k$ is odd.*

Turning to the symplectic case, let $V$ be a $2n$-dimensional complex vector space with a skew-symmetric bilinear form $a$. Then there exists a basis $e_i\ (i=1,\ldots,2n)$ of $V$ such that $a(e_i, e_{n+i})=1$ [hence $a(e_{n+i}, e_i)=-1$] for $i=1,\ldots,n$ and $a(e_i, e_j)=0$ if $|i-j|\neq n$. Introducing coordinates for $g\in\mathrm{Sp}_{2n}$ via $e_i g=\sum_{l=1}^{2n} g_{li}e_l$, we obtain $a(e_i, e_j g)=g_{i+n,j}$, $a(e_i, e_{j+n}g)=g_{i+n,j+n}$, $a(e_{i+n}, e_j g)=-g_{ij}$, $a(e_{i+n}, e_{j+n}g)=-g_{i,j+n}$. Given functions $\varphi,\ \psi\in\mathcal{F}(k,n)$, $\alpha,\ \beta\in\mathcal{F}(k,\{0,1\})$ we can write for $j=1,\ldots,k$

$$g_{j(\varphi+n\alpha),j(\psi+n\beta)}=(-1)^{1-j\alpha}a\big(e_{j(\varphi+n(1-\alpha))}, e_{j(\psi+n\beta)}g\big),$$

hence

(9)
$$\int \prod_{j=1}^k g_{j(\varphi+n\alpha),j(\psi+n\beta)}\,dg$$

$$= (-1)^{k-\sum_j j\alpha}\int a\big(e_{\varphi+n(1-\alpha)}, e_{\psi+n\beta}\rho_k(g)\big)\,dg$$

$$= (-1)^{k-\sum_j j\alpha}a\Big(e_{\varphi+n(1-\alpha)}, \int e_{\psi+n\beta}\rho_k(g)\,dg\Big),$$

where the extension of $a$ to the tensor power is defined as in (7). Since there exist no nontrivial invariants in odd tensor powers of $V$, we may assume that $k=2l$ is even.

To describe a basis of $[V^{\otimes 2l}]^{\mathrm{Sp}_{2n}}$, we have to proceed slightly more carefully than in the orthogonal case. An ordered pair partition of $\{1,\ldots,2l\}$ is a set $\mathfrak{m}^\mathfrak{o} = \{(m_\nu, n_\nu) : \nu = 1, \ldots, l\}$ of ordered pairs such that $\mathfrak{m} := \{\{m_\nu, n_\nu\} : \nu = 1, \ldots, l\} \in \mathcal{M}(2l)$. Set

$$\theta_{\mathfrak{m}^\mathfrak{o}} = \sum_{\eta \in \mathcal{F}(l,n)} \sum_{\varepsilon \in \mathcal{F}(l,\{0,1\})} \bigotimes_{j=1}^{2l} v(\eta, \varepsilon, j),$$

where

$$(10) \qquad v(\eta, \varepsilon, j) := \begin{cases} e_{\nu(\eta + n\varepsilon)}, & \text{if } j = m_\nu, \\ e_{\nu(\eta + n(1-\varepsilon))}(-1)^{\nu\varepsilon}, & \text{if } j = n_\nu. \end{cases}$$

For $n \geq l$, it follows from Weyl's First Fundamental Theorem and Section 3.4 of [16] that

$$\{\theta_{\mathfrak{m}^\mathfrak{o}} : \mathfrak{m}^\mathfrak{o} \in \mathcal{M}(2l)\}$$

is a basis for $[V^{\otimes 2l}]^{\mathrm{Sp}_{2n}}$, where $\mathcal{M}(2l)$ is identified with the set of all ordered pair partitions $\{(m_\nu, n_\nu) : \nu = 1, \ldots, l\}$ such that $m_\nu < n_\nu$ for all $\nu$.

Now let $\varphi \in \mathcal{F}(2l, n)$ and $\alpha \in \mathcal{F}(2l, \{0, 1\})$. We say that $\varphi + n\alpha \in \mathcal{F}(\mathfrak{m}^\mathfrak{o}, 2n)$ if $\varphi \in \mathcal{F}(\mathfrak{m}, n)$ and $\{m_\nu\alpha, n_\nu\alpha\} = \{0, 1\}$ for all $\nu = 1, \ldots, l$. Obviously, for $a(\theta_{\mathfrak{m}^\mathfrak{o}}, e_{\varphi + n\alpha})$ not to vanish it is necessary and sufficient that $\varphi + n\alpha \in \mathcal{F}(\mathfrak{m}^\mathfrak{o}, 2n)$. In this case, there exist unique $\eta, \varepsilon$ such that

$$a\left(\bigotimes_{j=1}^{2l} v(\eta, \varepsilon, j), e_{\varphi + n\alpha}\right) \neq 0.$$

Specifically, $\eta = \varphi$ and $\nu\varepsilon = 1$ if, and only if, $m_\nu\alpha = 0$. For these $\eta, \varepsilon$ we have

$$(11) \qquad a(\theta_{\mathfrak{m}^\mathfrak{o}}, e_{\varphi + n\alpha}) = a\left(\bigotimes_{j=1}^{2l} v(\eta, \varepsilon, j), e_{\varphi + n\alpha}\right) = \pm 1.$$

To determine the correct sign, let us disregard the sign in (10) for a moment. Any $\nu = 1, \ldots, l$ yields a factor of the form $a(e_{i+n}, e_i)a(e_i, e_{i+n}) = -1$ in (11). So we have a factor $(-1)^l$ regardless of the specific choice for $\varepsilon$. Now, we get an extra minus sign from (10) if $\nu\varepsilon = 1$, hence $m_\nu\alpha = 0$. So we see that

$$(12) \qquad (11) = (-1)^{l + \#\{\nu : m_\nu\alpha = 0\}}.$$

It remains to compute (9). To this end, observe that $k = 2l$ being even, $a$ is a symmetric bilinear form on $V^{\otimes k}$. So the corresponding steps of the orthogonal case are available. Applying (12) to $\varphi + n(1 - \alpha)$ and $\psi + n\beta$ and keeping in mind that $k = 2l$ is even, we see that (9) equals

$$(13) \qquad (-1)^{\#\{j : j\alpha = 1\}}$$
$$\times \sum_{\mathfrak{m}^\mathfrak{o}, \mathfrak{n}^\mathfrak{o} \in \mathcal{M}(2l)} \mathrm{Wg}_{\mathrm{Sp}_{2n}}(\mathfrak{m}^\mathfrak{o}, \mathfrak{n}^\mathfrak{o})(-1)^{\#\{\nu : m_\nu\alpha = 1\} + \#\{\nu : m'_\nu\beta = 0\}}$$
$$\times 1_{\mathcal{F}(\mathfrak{m}^\mathfrak{o}, 2n)}(\varphi + n\alpha) 1_{\mathcal{F}(\mathfrak{n}^\mathfrak{o}, 2n)}(\psi + n\beta),$$

where we have generically written

$$(14) \qquad \mathfrak{m}^{\mathfrak{o}} = \{(m_\nu, n_\nu) : \nu = 1, \dots, l\}, \qquad \mathfrak{n}^{\mathfrak{o}} = \{(m'_\nu, n'_\nu) : \nu = 1, \dots, l\}.$$

Observe that all summands vanish unless the prefactor (13) equals $(-1)^l$. So the overall sign factor is

$$(-1)^{l + \#\{\nu : m_\nu \alpha = 1\} + \#\{\nu : m'_\nu \beta = 0\}} = (-1)^{\#\{\nu : m_\nu \alpha = 0\} + \#\{\nu : m'_\nu \beta = 0\}}.$$

Summing up, we have the following:

PROPOSITION 3.2. *Let $\varphi, \psi \in \mathcal{F}(k, n)$, $\alpha, \beta \in \mathcal{F}(k, \{0, 1\})$. Then*

$$\int_{\mathrm{Sp}_{2n}} \prod_{j=1}^{k} g_{j(\varphi + n\alpha), j(\psi + n\beta)} \lambda_{\mathrm{Sp}_{2n}}(dg)$$

$$= \sum_{\mathfrak{m}^{\mathfrak{o}}, \mathfrak{n}^{\mathfrak{o}} \in \mathcal{M}(k)} \mathrm{Wg}_{\mathrm{Sp}_{2n}}(\mathfrak{m}^{\mathfrak{o}}, \mathfrak{n}^{\mathfrak{o}})(-1)^{\#\{\nu : m_\nu \alpha = 0\} + \#\{\nu : m'_\nu \beta = 0\}}$$

$$\times 1_{\mathcal{F}(\mathfrak{m}^{\mathfrak{o}}, 2n)}(\varphi + n\alpha) 1_{\mathcal{F}(\mathfrak{n}^{\mathfrak{o}}, 2n)}(\psi + n\beta),$$

*where we adopt the notation of (14) above and identify $\mathcal{M}(k)$ $(k = 2l)$ with the set of all ordered pair partitions $\{(m_\nu, n_\nu) : \nu = 1, \dots, l\}$ such that $m_\nu < n_\nu$ for all $\nu$.*

REMARK 3.3. Since $gJ = J\overline{g}$ for $g \in \mathrm{Sp}_{2n}$, Proposition 3.2 covers also the case of complex conjugates in the integrand.

Finally, let us record the unitary integration formula.

PROPOSITION 3.4. *Let $\varphi, \psi, \varphi', \psi' \in \mathcal{F}(k, n)$. Then*

$$\int_{\mathrm{U}_n} \prod_{j=1}^{k} u_{j\varphi, j\psi} \overline{u_{j\varphi', j\psi'}} \; \lambda_{\mathrm{U}_n}(du) = \sum_{\sigma, \tau \in \mathrm{S}_k} \delta(\varphi, \sigma\varphi') \delta(\psi, \tau\psi') \mathrm{Wg}_{\mathrm{U}_n}(\tau\sigma^{-1}),$$

*where the unitary Weingarten function $\mathrm{Wg}_{\mathrm{U}_n}$ is the inverse for convolution on the symmetric group $\mathrm{S}_k$ of the function $\sigma \mapsto n^{\#\mathrm{cycles}(\sigma)}$.*

In what follows, there will be no need to work with the definitions of the various Weingarten functions. It suffices to know their fundamental asymptotic properties:

PROPOSITION 3.5.

$$(15) \qquad \mathrm{Wg}_{\mathrm{U}_n}(\sigma) = \begin{cases} n^{-k}(1 + O(n^{-1})), & \text{if } \sigma = \mathrm{id}_k, \\ O(n^{-(k+1)}), & \text{otherwise.} \end{cases}$$

*For $k = 2l$, we have*

$$
(16) \qquad \mathrm{Wg}_{\mathrm{O}_n}(\mathfrak{m}_1, \mathfrak{m}_2) = \begin{cases} n^{-l}(1 + O(n^{-1})), & \text{if } \mathfrak{m}_1 = \mathfrak{m}_2 \in \mathcal{M}(k), \\ O(n^{-(l+1)}), & \text{otherwise,} \end{cases}
$$

$$
(17) \qquad \mathrm{Wg}_{\mathrm{Sp}_{2n}}(\mathfrak{m}_1^{\mathfrak{o}}, \mathfrak{m}_2^{\mathfrak{o}}) = \begin{cases} n^{-l}(1 + O(n^{-1})), & \text{if } \mathfrak{m}_1 = \mathfrak{m}_2 \in \mathcal{M}(k), \\ O(n^{-(l+1)}), & \text{otherwise.} \end{cases}
$$

PROOF. [4], Theorem 3.13. □

**4. The proofs.** Recall that we write $W_{\mathcal{C}}^{(n)}$ for the $*$-closed $\mathbb{K}$-linear hull of $\mathcal{S}_{\mathcal{C}}^{(n)}$ in $\mathbb{K}^{\delta n \times \delta n}$, and $P_{\mathcal{C}}^{(n)}$ for the orthogonal projection of $\mathbb{K}^{\delta n \times \delta n}$ onto $W_{\mathcal{C}}^{(n)}$ with respect to the scalar product $(A, B) \mapsto \mathrm{Tr}(AB^*)$. Then obviously $\mathrm{Re}(\mathrm{Tr}(A^{(\nu,n)} V_{\mathcal{C}}^{(n)})) = \mathrm{Re}(\mathrm{Tr}(P_{\mathcal{C}}^{(n)}(A^{(\nu,n)}) V_{\mathcal{C}}^{(n)}))$, so we can restrict our attention to $A^{(\nu,n)}$ taken from $W_{\mathcal{C}}^{(n)}$. Observe that the map $W_{\mathcal{C}}^{(n)} \to \mathbb{R}, A^{(\nu,n)} \mapsto \mathrm{Re}(\mathrm{Tr}(A^{(\nu,n)} V_{\mathcal{C}}^{(n)}))$ is $\mathbb{R}$-linear. By Cramér–Wold (see, e.g., [2], Theorem 29.4), this implies that it suffices to prove convergence to a Gaussian limit of variance $\lim_n \frac{\gamma_{\mathcal{C}}}{n} \mathrm{Tr}(A^{(\nu,n)}(A^{(\nu,n)})^*) = \lim_n \frac{\gamma_{\mathcal{C}}}{n} \mathrm{Tr}[P_{\mathcal{C}}^{(n)}(A^{(\nu,n)})(P_{\mathcal{C}}^{(n)}(A^{(\nu,n)}))^*]$ for an individual $\nu$. Recall that for a Gaussian distribution, the third and higher cumulants vanish. Denote the $r$th cumulant of $\mathrm{Re}(\mathrm{Tr}(A^{(\nu,n)} V_{\mathcal{C}}^{(n)}))$ by $C_r^{(n)}$. It follows from [4], Theorem 3.16, that for $r \geq 3$, $C_r^{(n)} = o(C_2^{(n)})$, as $n \to \infty$. So it suffices to show that $C_2^{(n)}$, that is, the variance, converges to the correct finite limit. We verify this in a case-by-case fashion.

**4.1. Class A.** We have $\Phi_{\mathrm{A}}^{(n)}(g) = g (g \in \mathrm{U}_n)$ and $W_{\mathrm{A}}^{(n)} = \mathbb{C}^{n \times n}$. By abuse of language, we also write $g = g^{(n)}$ for a random element of $\mathrm{U}_n$, with distribution given by Haar measure. Haar measure on $\mathrm{U}_n$ being invariant under multiplication by scalars $z \in \mathbb{C}$ with $|z| = 1$, we see that $\mathbb{E}(\mathrm{Re}\,\mathrm{Tr}(A^{(n)} g^{(n)})) = 0$ for all $A^{(n)} \in W_{\mathrm{A}}^{(n)}$ $(n \in \mathbb{N})$. Now,

$$
(18) \qquad \begin{aligned} \mathbb{E}((\mathrm{Re}\,\mathrm{Tr}(A^{(n)} g^{(n)}))^2) &= \mathbb{E}\left( \left( \frac{\mathrm{Tr}(A^{(n)} g^{(n)}) + \overline{\mathrm{Tr}(A^{(n)} g^{(n)})}}{2} \right)^2 \right) \\ &= \frac{1}{2} \mathbb{E}(\mathrm{Tr}(A^{(n)} g^{(n)}) \overline{\mathrm{Tr}(A^{(n)} g^{(n)})}). \end{aligned}
$$

The second equality is evident from Proposition 3.4. Now, (18) equals

$$
\frac{1}{2} \sum_{i,j,i',j'=1}^n a_{ij}^{(n)} \overline{a_{i'j'}^{(n)}} \int g_{ji} \overline{g_{j'i'}} \; \lambda_{\mathrm{U}_n}(dg) = \frac{1}{2} \mathrm{Wg}_{\mathrm{U}_n}(\mathrm{id}_1) \sum_{i,j=1}^n a_{ij}^{(n)} \overline{a_{ij}^{(n)}}
$$

by Proposition 3.4. By Proposition 3.5, this is asymptotically equivalent to

$$
\frac{1}{2n} \mathrm{Tr}(A^{(n)}(A^{(n)})^*).
$$

4.2. *Class AI.* To simplify notation, we will from now on suppress dependence on $n$ in most cases. We have $\Phi_{\mathrm{AI}}(g) = gg' (g \in \mathrm{U}_n), W_{\mathrm{AI}} = \{A \in \mathbb{C}^{n \times n} : A' = A\}$. Recall from Section 2 that we write $V = V^{(n)} = V_{\mathrm{AI}}^{(n)}$ for a random element of $\mathcal{S}_{\mathrm{AI}}^{(n)}$, distributed according to the push-forward of Haar measure under $\Phi_{\mathrm{AI}}^{(n)}$. Continuing to abuse language, and to write $g$ for a Haar distributed element of $\mathrm{U}_n$, we have $V = gg'$. As in the A case, one sees that $\mathbb{E}(\mathrm{ReTr}(AV)) = 0$ and

$$\mathbb{E}((\mathrm{ReTr}(AV))^2) = \tfrac{1}{2}\mathbb{E}(\mathrm{Tr}(AV)\overline{\mathrm{Tr}(AV)}).$$

Now this equals

$$\tfrac{1}{2} \sum_{i,j,k,i',j',k'} a_{ji}\overline{a_{j'i'}} \int g_{ik}g_{jk}\overline{g_{i'k'}}\,\overline{g_{j'k'}}\,dg,$$

which is asymptotically equivalent to

$$\mathrm{Wg}_{\mathrm{U}_n}(\mathrm{id}_2)\frac{n}{2}\left(\sum_{i,j} a_{ji}\overline{a_{ji}} + \sum_{i,j} a_{ji}\overline{a_{ij}}\right),$$

hence to

$$\frac{1}{n}\mathrm{Tr}(A^{(n)}(A^{(n)})^*),$$

since $A' = A$.

4.3. *Class AII.* Here $\Phi_{\mathrm{AII}}(g) = gJ'g'J$ ($g \in \mathrm{U}_{2n}$) and $W_{\mathrm{AII}} = \{A \in \mathbb{C}^{2n \times 2n} : (AJ)' = -AJ\}$. As in the A and AI cases, $\mathbb{E}(\mathrm{ReTr}(A^{(n)}V^{(n)})) = 0$ and

$$\mathbb{E}((\mathrm{ReTr}(A^{(n)}V^{(n)}))^2) = \tfrac{1}{2}\mathbb{E}(\mathrm{Tr}(A^{(n)}V^{(n)})\overline{\mathrm{Tr}(A^{(n)}V^{(n)})}).$$

Write $J = (J_{ij})$. One has $\mathrm{Tr}(AV) = \mathrm{Tr}(AgJ'g'J) = -\mathrm{Tr}(BgJg')$ with $B := JA$. Now,

$$\tfrac{1}{2}\int \mathrm{Tr}(BgJg')\overline{\mathrm{Tr}(BgJg')}\,dg$$

$$= \tfrac{1}{2}\sum_{i,j,k,l,i',j',k',l'} b_{ij}J_{kl}\overline{b_{i'j'}}J_{k'l'}\int g_{jk}g_{il}\overline{g_{j'k'}}\,\overline{g_{i'l'}}\,dg.$$

This is asymptotically equivalent to

$$\tfrac{1}{2}\mathrm{Wg}_{\mathrm{U}_{2n}}(\mathrm{id}_2)\left(\sum_{i,j,k,l} b_{ij}J_{kl}\overline{b_{ij}}J_{kl} + \sum_{i,j,k,l} b_{ij}J_{kl}\overline{b_{ji}}J_{lk}\right),$$

hence to

$$\frac{1}{2(2n)^2}(\mathrm{Tr}(BB^*)\,\mathrm{Tr}(JJ') + \mathrm{Tr}(B\overline{B})\,\mathrm{Tr}(J^2)) = \frac{1}{(2n)^2}\mathrm{Tr}(AA^*)\,\mathrm{Tr}(JJ').$$

So we have obtained that $\mathbb{E}((\mathrm{ReTr}(A^{(n)}V^{(n)}))^2) \sim \frac{1}{2n}\mathrm{Tr}(A^{(n)}(A^{(n)})^*)$.

4.4. *Class AIII.* In this case one has $\Phi_{\text{AIII}}(g) = g I_{p,q} g^* I_{p,q}, g \in U_n$ ($n = p + q$) and $W_{\text{AIII}} = \{A \in \mathbb{C}^{n \times n} : A I_{p,q} = (A I_{p,q})^*\}$. For $A = A^{(n)} \in W_{\mathcal{C}}^{(n)}$, we have $\operatorname{Tr}(A V^{(n)}) = \operatorname{Tr}(A g I_{p,q} g^* I_{p,q}) = \operatorname{Tr}(B g I_{p,q} g^*)$ with $B := I_{p,q} A$. Observe that this is a real random variable. Writing $I_{p,q} = (I_{ij})$, we have that $\mathbb{E}(\operatorname{Tr}(B g I_{p,q} g^*)) = \sum_{i,j,k,l} b_{ij} I_{kl} \int g_{jk} \overline{g_{il}}\, dg$ is asymptotically equivalent to

$$\frac{1}{n} \sum_{i,k} b_{ii} I_{kk} = \frac{1}{n} \operatorname{Tr}(B) \operatorname{Tr}(I_{p,q})$$

(19)

$$= \frac{p(n) - q(n)}{n} \operatorname{Tr}(B) = \frac{p(n) - q(n)}{n} \operatorname{Tr}(A^{(n)} I_{p(n),q(n)}).$$

On the other hand,

$$\mathbb{E}(\operatorname{Tr}(B g I_{p,q} g^*))^2 = \sum_{i,j,k,l,i',j',k',l'} b_{ij} b_{i'j'} I_{kl} I_{k'l'} \int g_{jk} g_{j'k'} \overline{g_{il}}\,\overline{g_{i'l'}}\, dg$$

is asymptotically equivalent to

$$\operatorname{Wg}_{U_n}(\operatorname{id}_2)\left( \sum_{i,i',k,k'} b_{ii} b_{i'i'} I_{kk} I_{k'k'} + \sum_{i,i',l,l'} b_{ii'} b_{i'i} I_{l'l} I_{ll'} \right),$$

hence to

$$\frac{1}{n^2}((\operatorname{Tr} B)^2 (\operatorname{Tr} I_{p,q})^2 + \operatorname{Tr}(B^2) \operatorname{Tr}(I_{p,q}^2)).$$

Now (19), the assumptions of Theorem 2, and the definitions of $B$ and $W_{\text{AIII}}$ imply that as $n \to \infty$,

$$\mathbb{V}(\operatorname{Tr}(B^{(n)} g^{(n)} I_{p,q}(g^{(n)})^*)) = \mathbb{V}(\operatorname{Tr}(A^{(n)} V^{(n)})) \sim \frac{1}{n} \operatorname{Tr}(A^{(n)}(A^{(n)})^*).$$

4.5. *Class B/D.* Here $\Phi_{\text{B/D}}(g) = g, g \in O_n$, and $W_{\text{B/D}} = \mathbb{R}^{n \times n}$. By invariance of Haar measure for $A \in W_{\text{B/D}}$ one has $\mathbb{E}(\operatorname{Tr}(Ag)) = 0$. On the other hand, by Proposition 3.1

$$\mathbb{E}(\operatorname{Tr}(Ag)^2) = \sum_{i,j,k,l} a_{ij} a_{kl} \int g_{ji} g_{lk}\, dg = \operatorname{Wg}_{O_n}(\{\{1,2\}\}, \{\{1,2\}\}) \sum_{i,j} a_{ij}^2,$$

which is by Proposition 3.5 asymptotically equivalent to $\frac{1}{n} \operatorname{Tr}(A^{(n)}(A^{(n)})')$.

4.6. *Class BDI.* In this case one has $\Phi_{\text{BDI}}(g) = g I_{p,q} g' I_{p,q}$, $g \in O_n$, $n = p + q$ and $W_{\text{BDI}} = \{A \in \mathbb{R}^{n \times n} : (A I_{p,q})' = A I_{p,q}\}$. Then $\operatorname{Tr}(A V^{(n)}) = \operatorname{Tr}(A g I_{p,q} g' I_{p,q}) = \operatorname{Tr}(B g I_{p,q} g')$ with $B := I_{p,q} A$. Writing again $I_{p,q} = (I_{ij})$, one obtains

$$\mathbb{E}(\operatorname{Tr}(B g I_{p,q} g')) = \sum_{i,j,k,l} b_{ij} I_{kl} \int g_{jk} g_{il}\, dg,$$

which is asymptotically equivalent to

$$(20) \quad \frac{1}{n} \sum_{i,k} b_{ii} I_{kk} = \frac{1}{n} \operatorname{Tr}(B) \operatorname{Tr}(I_{p,q}) = \frac{p(n) - q(n)}{n} \operatorname{Tr}(A^{(n)} I_{p(n),q(n)}).$$

On the other hand,

$$\mathbb{E}((\operatorname{Tr}(BgI_{p,q}g'))^2) = \sum_{i,j,k,l,i',j',k',l'} b_{ij} b_{i'j'} I_{kl} I_{k'l'} \int g_{jk} g_{il} g_{j'k'} g_{i'l'} \, dg.$$

By Propositions 3.1 and 3.5, this is asymptotically equivalent to

$$(21) \quad \frac{1}{n^2} \left( \sum_{i,i',k,k'} b_{ii} b_{i'i'} I_{kk} I_{k'k'} + \sum_{i,i',k,k'} b_{ii'} b_{i'i} I_{kk'} I_{k'k} + \sum_{i,j,k,l} b_{ij} b_{ij} I_{kl} I_{kl} \right)$$

$$= \frac{1}{n^2} ((\operatorname{Tr} B)^2 (\operatorname{Tr} I_{p,q})^2 + \operatorname{Tr}(B^2) \operatorname{Tr}(I_{p,q}^2) + \operatorname{Tr}(BB') \operatorname{Tr}(I_{p,q} I'_{p,q}))$$

$$= \frac{1}{n^2} ((\operatorname{Tr} B)^2 (p - q)^2 + 2n \operatorname{Tr}(BB')),$$

where the three sums in (21) correspond to the pair partitions $\{\{1,2\}, \{3,4\}\}$, $\{\{1,4\}, \{2,3\}\}$ and $\{\{1,3\}, \{2,4\}\}$ of $\{1,2,3,4\}$. Comparing with (20), we see that $\mathbb{V}(\operatorname{Tr}(A^{(n)} V^{(n)})) \sim \frac{2}{n} \operatorname{Tr}(A^{(n)} (A^{(n)})')$.

4.7. *Class DIII.* Here one has $\Phi_{\text{DIII}}(g) = gJ'g'J, g \in O_{2n}$, and $W_{\text{DIII}} = \{A \in \mathbb{R}^{2n \times 2n} : (AJ)' = -AJ\}$. Then $\operatorname{Tr}(A^{(n)} V^{(n)}) = \operatorname{Tr}(AgJ'g'J) = \operatorname{Tr}(BgJ'g')$ with $B := JA$.

$$\mathbb{E}(\operatorname{Tr}(BgJ'g')) = \sum_{i,j,k,l} b_{ij} J_{lk} \int g_{jk} g_{il} \, dg = \frac{1}{n} \sum_{i,k} b_{ii} J_{kk} = 0$$

by definition of $J$. On the other hand,

$$\mathbb{E}((\operatorname{Tr}(BgJ'g'))^2) = \sum_{i,j,k,l,i',j',k',l'} b_{ij} J_{lk} b_{i'j'} J_{l'k'} \int g_{jk} g_{il} g_{j'k'} g_{i'l'} \, dg.$$

As in the BDI case, one is led to sums corresponding to the pair partitions of $\{1,2,3,4\}$. The sum corresponding to $\{\{1,2\}, \{3,4\}\}$ vanishes, because it involves diagonal terms of $J$. The sum corresponding to $\{\{1,3\}, \{2,4\}\}$ is $\frac{1}{(2n)^2} \sum_{i,j,k,l} b_{ij} b_{ij} J_{lk} J_{lk} = \frac{1}{(2n)^2} \operatorname{Tr}(BB') \operatorname{Tr}(JJ') = \frac{1}{2n} \operatorname{Tr}(AA')$. The last sum is $\frac{1}{(2n)^2} \sum_{i,i',l,l'} b_{ii'} J_{ll'} b_{i'i} J_{l'l} = \frac{1}{(2n)^2} \operatorname{Tr}(B^2) \operatorname{Tr}(J^2) = (-\operatorname{Tr}(AA'))(-\frac{1}{2n})$. Summing up,

$$\mathbb{E}((\operatorname{Tr}(A^{(n)} V^{(n)}))^2) \sim \frac{1}{n} \operatorname{Tr}(A^{(n)} (A^{(n)})').$$

4.8. *Class C.* Here we have $\Phi_{\mathrm{C}}(g) = g, g \in \mathrm{Sp}_{2n}$, and $W_{\mathrm{C}} = \mathbb{H}^{n \times n}$. Since $\mathbb{H}^{n \times n}$ is closed under multiplication, we see that $\mathrm{Tr}(Ag)$ is real for

$$(22) \qquad A = \begin{pmatrix} X & Y \\ -\overline{Y} & \overline{X} \end{pmatrix} \in \mathbb{H}^{n \times n},$$

and the same argument applies to classes CI and CII below. Invariance of Haar measure yields $\mathbb{E}(\mathrm{Tr}(Ag)) = 0$. Furthermore, from Proposition 3.2 we obtain

$$\mathbb{E}((\mathrm{Tr}(A^{(n)}V^{(n)}))^2)$$

$$= \mathbb{E}((\mathrm{Tr}(A^{(n)}g))^2)$$

$$= \sum_{\substack{i,j,i',j' \\ \in \{1,\dots,n\}}} \sum_{s,t,s',t' \in \{0,1\}} a_{j+nt,i+ns} a_{j'+nt',i'+ns'} \int g_{i+ns,j+nt} g_{i'+ns',j'+nt'} \, dg$$

$$= \mathrm{Wg}_{\mathrm{Sp}_{2n}}(\{(1,2)\},\{(1,2)\}) \sum_{s,t} (-1)^{s+t} \sum_{i,j} a_{j+nt,i+ns} a_{j+n(1-t),i+n(1-s)}$$

$$\sim \frac{1}{n} \sum_{i,j} a_{ji} a_{j+n,i+n} + a_{j+n,i+n} a_{ji} - a_{j+n,i} a_{j,i+n} - a_{j,i+n} a_{j+n,i}$$

$$= \frac{2}{n}(\mathrm{Tr}(X(\overline{X})') - \mathrm{Tr}(Y(-\overline{Y})')) = \frac{2}{n} \mathrm{Tr}(XX^* + YY^*)$$

$$= \frac{1}{n} \mathrm{Tr}(A^{(n)}(A^{(n)})^*).$$

4.9. *Class CI.* We have $\Phi_{\mathrm{CI}}(g) = gI_{n,n}g^{-1}I_{n,n} (g \in \mathrm{Sp}_{2n})$ and $W_{\mathrm{CI}} = \{A \in \mathbb{H}^{n \times n} : (AI_{n,n})^* = AI_{n,n}\}$. Set

$$K := JI_{n,n} = -I_{n,n}J = \begin{pmatrix} 0 & I_n \\ I_n & 0 \end{pmatrix}$$

and write $K = (\kappa_{ij})$. Since $gJ = J\overline{g}$ for $g \in \mathrm{Sp}_{2n}$, we have for any $A \in W_{\mathrm{CI}}$ that

$$\mathrm{Tr}(A^{(n)}V^{(n)})$$

$$= \mathrm{Tr}(AgI_{n,n}g^{-1}I_{n,n})$$

$$= \mathrm{Tr}(AgI_{n,n}(-J)g'JI_{n,n})$$

$$= \mathrm{Tr}(AgKg'K) = \mathrm{Tr}(BgKg')$$

$$= \sum_{i,j,k,l \in \{1,\dots,n\}} \sum_{q,r,s,t \in \{0,1\}} b_{i+nq,j+nr} g_{j+nr,k+ns} \kappa_{k+ns,l+nt} g_{i+nq,l+nt},$$

where $KA =: B = (b_{ij})$. It follows easily from Proposition 3.2 and the symmetry of $K$ that $\mathbb{E}(\mathrm{Tr}(BgKg')) = 0$. Turning to $\mathbb{E}((\mathrm{Tr}(BgKg'))^2)$, we apply

Proposition 3.2 to

$$\sum_{\substack{i,j,k,l \\ i',j',k',l'}} \sum_{\substack{q,r,s,t \\ q',r',s',t'}} b_{i+nq,j+nr} b_{i'+nq',j'+nr'} \kappa_{k+ns,l+nt} \kappa_{k'+ns',l'+nt'}$$

$$\times \int g_{j+nr,k+ns} g_{i+nq,l+nt} g_{j'+nr',k'+ns'} g_{i'+nq',l'+nt'} \, dg.$$

The ordered pair partition $\{(1,2),(3,4)\}$ leads to summands of the form

$$(-1)^{r+s+r'+s'} b_{i+n(1-r),i+nr} b_{i'+n(1-r'),i'+nr'} \kappa_{k+ns,k+n(1-s)} \kappa_{k'+ns',k'+n(1-s')}.$$

$K$ being symmetric, the summands come in $\pm$ pairs, and the total sum vanishes. To $\{(1,3),(2,4)\}$ correspond summands of the form

$$(-1)^{q+r+s+t} b_{i+nq,j+nr} b_{i+n(1-q),j+n(1-r)} \kappa_{k+ns,l+nt} \kappa_{k+n(1-s),l+n(1-t)}.$$

In view of the specific form of $K$, for such a summand not to vanish it is necessary that $k = l$ and $s \neq t$. Hence the contribution from $\{(1,3),(2,4)\}$ is

$$(23) \qquad \begin{aligned} & 2n \frac{1}{n^2} \sum_{i,j,q,r} (-1)^{q+r+1} b_{i+nq,j+nr} b_{i+n(1-q),j+n(1-r)} \\ & = \frac{2}{n} \sum_{i,j} -b_{ij} b_{i+n,j+n} + b_{i,j+n} b_{i+n,j} + b_{i+n,j} b_{i,j+n} - b_{i+n,j+n} b_{i,j}. \end{aligned}$$

Note that we have obtained the factor $2n$ since we have suppressed the sums over $k$ and $s$. Writing $A$ as in (22) above, we have

$$B = \begin{pmatrix} -\overline{Y} & \overline{X} \\ X & Y \end{pmatrix},$$

and since $AI_{n,n}$ was assumed to be hermitian, we have $\overline{X}' = X$ and $Y = Y'$. This means $b_{i+n,j+n} = b_{j+n,i+n} = -\overline{b_{ij}}$ and $b_{i,j+n} = \overline{b_{i+n,j}} = \overline{b_{j,i+n}}$. Hence, (23) equals

$$\frac{4}{n}(\operatorname{Tr}(XX^*) + \operatorname{Tr}(YY^*)) = \frac{2}{n} \operatorname{Tr}(AA^*).$$

Finally, the ordered pair partition $\{(1,4),(2,3)\}$ yields summands of the form

$$(-1)^{r+s+q+t} b_{i+nq,j+nr} b_{j+n(1-r),i+n(1-q)} \kappa_{k+ns,l+nt} \kappa_{l+n(1-t),k+n(1-s)},$$

and an analogous argument shows that the contribution of $\{(1,4),(2,3)\}$ is

$$\frac{2}{n} \sum_{ij} b_{i,j+n} b_{j,i+n} + b_{i+n,j} b_{j+n,i+n} - b_{i,j} b_{j+n,i+n} - b_{i+n,j+n} b_{j,i} = \frac{2}{n} \operatorname{Tr}(AA^*).$$

Summing up, $\mathbb{E}((\operatorname{Tr}(A^{(n)} V^{(n)}))^2) = \mathbb{E}((\operatorname{Tr}(BgKg'))^2) \sim \frac{4}{n} \operatorname{Tr}(A^{(n)}(A^{(n)})^*).$

4.10. *Class CII.* We have $\Phi_{\text{CII}}(g) = g\mathcal{I}_{p,q}g^{-1}\mathcal{I}_{p,q}$ $(g \in \text{Sp}_{2n},\ n = p+q)$, $W_{\text{CII}} = \{A \in \mathbb{H}^{n \times n} : (A\mathcal{I}_{p,q})^* = A\mathcal{I}_{p,q}\}$. Write

$$K_{p,q} := J\mathcal{I}_{p,q} = \begin{pmatrix} 0 & -I_{p,q} \\ I_{p,q} & 0 \end{pmatrix} = \mathcal{I}_{p,q}J$$

and $K_{p,q} = (\kappa_{ij})$. Fix

$$A = \begin{pmatrix} X & Y \\ -\overline{Y} & \overline{X} \end{pmatrix} \in W_{\text{CII}}.$$

Then

$$\text{Tr}(AV) = \text{Tr}(A g\mathcal{I}_{p,q}g^{-1}\mathcal{I}_{p,q}) = \text{Tr}(A g\mathcal{I}_{p,q}(-J)g'J\mathcal{I}_{p,q}) = -\text{Tr}(BgK_{p,q}g'),$$

where $(b_{ij}) = B := K_{p,q}A$. Note that

$$(24) \qquad B = \begin{pmatrix} I_{p,q}\overline{Y} & -I_{p,q}\overline{X} \\ I_{p,q}X & I_{p,q}Y \end{pmatrix}.$$

Since $A\mathcal{I}_{p,q}$ is assumed to be Hermitian, we obtain the useful relations

$$(25) \qquad XI_{p,q} = (XI_{p,q})^* = I_{p,q}X^*, \qquad YI_{p,q} = -(YI_{p,q})' = -I_{p,q}Y'.$$

Now, $\mathbb{E}(-\text{Tr}(BgK_{p,q}g'))$ equals

$$-\sum_{i,j,k,l \in \{1,\ldots,n\}} \sum_{q,r,s,t \in \{0,1\}} b_{i+nq,j+nr}\kappa_{k+ns,l+nt} \int g_{j+nr,k+ns}g_{i+nq,l+nt}\,dg.$$

The leading order contribution comes from the case $j = i$, $k = l$, $q = (1-r)$, $t = (1-s)$, hence is asymptotically equivalent to

$$(26) \qquad -\frac{1}{n}\sum_{i,k}\sum_{r,s}(-1)^{r+s}b_{i+n(1-r),i+nr}\kappa_{k+ns,k+n(1-s)}$$

$$= -\frac{1}{n}\left(2\sum_{l=1}^{n}\kappa_{l,n+l}\right)(\text{Tr}(I_{p,q}(X+\overline{X})))$$

$$= 2\frac{(p(n)-q(n))}{n}\,\text{Tr}(A^{(n)}\mathcal{I}_{p(n),q(n)}).$$

Turning to the second moment, we have that $\mathbb{E}((\text{Tr}(BgK_{p,q}g'))^2)$ equals

$$\sum_{\substack{i,j,k,l \\ i',j',k',l'}} \sum_{\substack{q,r,s,t \\ q',r',s',t'}} b_{i+nq,j+nr}b_{i'+nq',j'+nr'}\kappa_{k+ns,l+nt}\kappa_{k'+ns',l'+nt'}$$

$$\times \int g_{j+nr,k+ns}g_{i+nq,l+nt}g_{j'+nr',k'+ns'}g_{i'+nq',l'+nt'}\,dg.$$

First of all, the leading order contribution corresponding to the ordered pair partition $\{(1,2),(3,4)\}$ is asymptotically equivalent to

$$\frac{1}{n^2} \sum_{i,k,i',k'} \sum_{r,s,r',s'} (-1)^{r+s+r'+s'} b_{i+n(1-r),i+nr} b_{i'+n(1-r'),i'+nr'}$$

$$\times \kappa_{k+ns,k+n(1-s)} \kappa_{k'+ns',k'+n(1-s')},$$

which coincides with the square of (26). So, asymptotically, the variance will come from the remaining ordered pair partitions. Next, from $\{(1,3),(2,4)\}$, we get summands of the form

$$(-1)^{r+s+q+t} b_{i+nq,j+nr} b_{i+n(1-q),j+n(1-r)} \kappa_{k+ns,l+nt} \kappa_{k+n(1-s),l+n(1-t)}.$$

For such a summand not to vanish, it is necessary that $k = l$ and $s \neq t$. Note that in this case the product of the $\kappa$ terms is $-1$. So we obtain a contribution

$$2n\frac{1}{n^2} \sum_{i,j} \sum_{q,r} (-1)^{q+r} b_{i+nq,j+nr} b_{i+n(1-q),j+n(1-r)}$$

$$= \frac{4}{n} \sum_{i,j} b_{i,j} b_{n+i,n+j} - b_{i,n+j} b_{n+i,j}$$

$$= \frac{4}{n} (\mathrm{Tr}((I_{p,q}\overline{Y})(I_{p,q}Y)') + \mathrm{Tr}((I_{p,q}\overline{X})(I_{p,q}X)'))$$

$$= \frac{4}{n} (\mathrm{Tr}(XX^*) + \mathrm{Tr}(YY^*)) = \frac{2}{n} \mathrm{Tr}(AA^*).$$

Finally, as to $\{(1,4),(2,3)\}$ one argues similarly to obtain a contribution

$$2n\frac{1}{n^2} (-2\,\mathrm{Tr}((I_{p,q}\overline{Y})(I_{p,q}Y))$$

$$+ \mathrm{Tr}((I_{p,q}\overline{X})^2) + \mathrm{Tr}((I_{p,q}X)^2)),$$

and one invokes (25) to see that this in fact equals $\frac{2}{n}\mathrm{Tr}(A^{(n)}(A^{(n)})^*)$.

CNRS, UMR 5208
Institut Camille Jordan
Université Lyon 1
43 boulevard du 11 novembre 1918
69622 Villeurbanne cedex
France
and
Département de mathématiques
  et de statistique
Université d'Ottawa
585 King Edward
Ottawa, Ontario
Canada K1N 6N5
E-mail: collins@math.univ-lyon1.fr

Fakultät für Mathematik
NA 4/32
Ruhr-Universität Bochum
D-44780 Bochum
Germany
E-mail: michael.stolz@ruhr-uni-bochum.de